\documentclass[11pt,twoside,a4paper]{article}
\usepackage{graphicx}
\usepackage{amsfonts}
\usepackage{bbm}
\usepackage[leqno]{amsmath}
\usepackage{mathrsfs}
\usepackage{fancyhdr}
 \usepackage{mathrsfs,amsmath,amssymb,amsthm}
\usepackage[dvips]{color}
 \usepackage{amssymb}
 \usepackage{amsbsy}
 \usepackage[dvips]{color}
\usepackage{titlesec}
\titlespacing*{\section}{0pt}{0.5\baselineskip}{0.5\baselineskip}
\titleformat*{\subsubsection}{\it}
\titlespacing*{\subsection}{0pt}{0.3\baselineskip}{0.3\baselineskip}
\titlespacing*{\subsubsection}{0pt}{0.3\baselineskip}{0.3\baselineskip}

\allowdisplaybreaks

 \theoremstyle{definition}
\theoremstyle{remark}  \newtheorem{remark}{\noindent\mbox{Remark}}
 \theoremstyle{plain}
 \theoremstyle{plain}\newtheorem{lemma}{\noindent\mbox{Lemma}}
\theoremstyle{plain} \newtheorem{theorem}{\noindent\mbox{Theorem}}
 \theoremstyle{plain}
 \theoremstyle{plain}
\theoremstyle{definition} 

 \def\proof{\noindent{\it Proof.~~}}
 \def\qed{\hfill$\Box$\medskip}
 \def\rto{\rightarrow\infty}
\def\z{\left}
\def\y{\right}
 \def\no{\nonumber}
 \def\mb{\mathbf}

\begin{document}

 \title{\bf{Asymptotics of entries of products of nonnegative 2-by-2 matrices}}                

\author{Hua-Ming \uppercase{Wang}\footnote{Email:hmking@ahnu.edu.cn; School of Mathematics and Statistics, Anhui Normal University, Wuhu, 241003, China }}

\date{}
\maketitle%

\vspace{-.7cm}

\begin{center}
\begin{minipage}[c]{12cm}
\begin{center}\textbf{Abstract}  \end{center}
Let $M$ and $M_n,n\ge1$ be nonnegative 2-by-2 matrices such that $\lim_{n\rto}M_n=M.$ It is usually  hard to estimate the entries of $M_{k+1}\cdots M_{k+n}$ which are useful in many applications.
 In this paper, under a mild condition, we show that up to a multiplication of some positive constants, entries of $M_{k+1}\cdots M_{k+n}$ are asymptotically the same as $\xi_{k+1}\cdots \xi_{k+n},$  the product of the tails of a continued fraction which is related to the matrices $M_k,k\ge1,$ as $n\rto.$

\vspace{0.2cm}

\vspace{0.2cm}

\end{minipage}
\end{center}

\section{Introduction}\label{intr}

In applications, it is usually important to evaluate the entries of the product of nonnegative matrices. For example, the mean matrix of the number of individuals alive at the $n$-th generation of a multitype branching process in random or varying environments is indeed the product nonnegative matrices \cite{jon,va,wy}, the escape probability of the random walk in random or varying environments from an interval can be written in terms of  products nonnegative matrices \cite{br,hs20}, and  Bernoulli convolutions and Gibbs properties of
linearly representable measures are also closely related to product of nonnegative matrices \cite{ot}. We mention that though the product of matrices has been extensively studied, see \cite{co,dl,fs,hj,235,jon,ot} and  references therein,  few computable estimate of the entries of product of matrices can be found in the literature.

In this paper, we consider the product of 2-by-2 nonnegative matrices. In \cite{hs20}, it has been shown that up to a multiplication of certain constants, entries of the product of nonnegative 2-by-2 matrices are asymptotically the same as the product of spectral radii of those matrices,  however some strong conditions which are far from natural are required. Instead of the spectral radii of the matrices used in \cite{hs20}, under mild condition, we will show in this paper that in some uniform way, entries of the product of nonnegative 2-by-2 matrices are asymptotically the same as
the products of the tails of a continued fraction which is related to those matrices.

 Throughout, we assume
\begin{align}\label{mg}M:=\left( \begin{array}{cc}
  a & b \\
  d &\theta \\
 \end{array}\right) \text{ and } M_k:=\left( \begin{array}{cc}
  a_k & b_k \\
  d_k &\theta_k \\
 \end{array}\right), k\ge1\end{align}
are nonnegative matrices and let \begin{align*}
  \varrho:=\frac{a+\theta+\sqrt{(a+\theta)^2+4(bd- a\theta)}}{2},\ \varrho_1:=\frac{a+\theta-\sqrt{(a+\theta)^2+4(bd-a\theta)}}{2}
\end{align*}
be the eigenvalues of $M.$

To introduce the main result, we need to introduce further some continued fractions  which are closely related to $M_k,k\ge1.$ For $n\ge1,$ let $\beta_n=\frac{b_{n+1}}{b_n(b_{n+1}d_{n+1}-a_{n+1}\theta_{n+1})},$  $\alpha_n=\frac{a_nb_{n+1}+b_n\theta_{n+1}}{b_n(b_{n+1}d_{n+1}-a_{n+1}\theta_{n+1})}$ and set
\begin{align}\label{xicn}
   &\xi_n:={\frac{\beta_{n}}{\alpha_{n}}}\begin{array}{c}
                                \\
                               +
                             \end{array}\frac{\beta_{n+1}}{\alpha_{n+1 }}\begin{array}{c}
                                \\
                               +
                             \end{array}\frac{\beta_{n+2}}{\alpha_{n+2}}\begin{array}{c}
                                \\
                               +\cdots
                             \end{array}.
\end{align}
Throughout the paper we use the notations  $\mb e_1=(1,0),$ $\mb e_2=(0,1)$ and   denote by $\mb v^t$ the transpose of a vector $\mb v.$  We also adopt the convention that empty sum equals 0 and empty product equals identity.  The theorem below provides an approach to evaluate the entries of the products of matrices.
\begin{theorem}\label{pml} Suppose that $\lim_{k\rto}M_k=M,$   $a+\theta\ne0,$ $b\ne0$ and $bd\ne a\theta.$ Then $\exists k_0>0$ such that for $k\ge k_0$ and $i,j\in\{1,2\},$ we have
\begin{align}\label{mx}
&\quad\ \frac{\mb e_iM_{k+1}\cdots M_{k+n}\mb e_j^t}{\xi_{k+1}^{-1}\cdots\xi_{k+n}^{-1}} \rightarrow \psi(i,j,k),
\end{align}
uniformly in $k$ as $n\rto,$ where $\psi(1,1,k)=\frac{\varrho-\theta}{\varrho-\varrho_1},$ $\psi(1,2,k)=\frac{b}{\varrho-\varrho_1},$
$\psi(2,1,k)=\frac{\varrho}{\varrho-\varrho_1}\z(\frac{\theta_{k+1}}{b_{k+1}}-\xi_{k+1}\frac{\mathrm{det}(M_{k+1})}{b_{k+1}}\y)$ and $\psi(2,2,k)=\frac{b}{\varrho-\varrho_1}\z(\frac{\theta_{k+1}}{b_{k+1}}-\xi_{k+1}\frac{\mathrm{det}(M_{k+1})}{b_{k+1}}\y).$
Furthermore, if we assume further $\varrho\ge1,$ then for $k\ge k_0$ and $i=1,$ $j\in\{1,2\},$ with the above $\psi(i,j,k),$ we have
\begin{align}
\label{ms}
    &\frac{\sum_{s=1}^{n+1}\mb e_iM_{k+s}\cdots M_{k+n}\mb e_j^t}{\sum_{s=1}^{n+1}\xi_{k+s}^{-1}\cdots \xi_{k+n}^{-1}}\rightarrow\psi(i,j,k),
  \end{align}
 uniformly in $k$ as $n\rto.$
\end{theorem}
\begin{remark} \textbf{(i)} We see from Theorem \ref{pml} that in order to evaluate the entries of $M_{k+1}\cdots M_{k+n},$ it is enough to know the asymptotics of $\xi_n$ which is   a limit periodic continued fraction and has been extensively studied, see \cite{lor,lw,w19} and references therein.
  \textbf{(ii)} The uniform convergence in Theorem \ref{pml} is crucial in some applications. For example, it is required when studying the number of cutpoints of (2,1) random walks in certain varying environments and the number of regeneration times of a two-type branching process in varying environments, see \cite{sw21, wt21}. \textbf{(iii)} It is shown in \cite{hs20} that $\mb e_i M_k\cdots M_n \mb e_j^t\sim c(k) \varrho(M_k)\cdots \varrho(M_n)$ as $n\rto,$ where $c(k)>0$ is a constant and  $\varrho(M_j)$ is the spectral radius of matrix $M_j.$  Based on \cite{hs20}, the asymptotics of the distribution of the extinction time of a two-type branching process was studied in \cite{wy}. However, in both \cite{hs20} and \cite{wy}, it is required that $M_k\rightarrow M$ in a strong manner so that some strong conditions are required. In view of \eqref{mx}, if we replace $\varrho(M_n)$ by $\xi_n,$ then the conditions in \cite{wy} can be weakened in a large extent.
\end{remark}

Let us explain the idea to prove Theorem \ref{pml}. It is difficult to work directly with the matrices $M_k, k\ge1.$ Instead, for $k\ge1$ we construct some matrices
\begin{align}\label{dta} A_k:=\left( \begin{array}{cc}
  \tilde a_k & \tilde b_k \\
  \tilde d_k &0 \\
 \end{array}\right) \text{ with }
\tilde a_k=a_k+\frac{b_k\theta_{k+1}}{ b_{k+1}}, \tilde b_k= b_k,\tilde  d_k=d_k-\frac{a_k\theta_k}{b_k}, \end{align} which is quasi-similar to  $M_k$ in the following sense
\begin{align}\label{am}
  A_k=\Lambda_k^{-1}M_k\Lambda_{k+1}\text{ where }
\Lambda_k=\left(
                 \begin{array}{cc}
                   1 & 0 \\
                   \theta_k/b_k & 1 \\
                 \end{array}
               \right).
\end{align}
In this way, it is observed that entries of $A_{k}\cdots A_n$ can be written in terms of  approximants of certain continued fractions, see \eqref{cpn} and Lemma \ref{axc} below. Then by analyzing subtly the limiting behaviours of the continued fraction and their tails and approximants, we can show that  for $i,j\in\{1,2\}$ and $k$ large enough,  there exists constant $ \phi(i,j,k)>0 $ such that
  $\frac{\mb e_iA_{k+1}\cdots A_{k+n}\mb e_j^t}{\xi_{k+1}^{-1}\cdots\xi_{k+n}^{-1}}\rightarrow \phi(i,j,k)
  $ uniformly in $k,$ see  Theorem \ref{s} below.  As a result, we can show the first part of Theorem \ref{pml} by using the connection of $M_k$ and $A_k$ given in \eqref{am}.

As far as the proof of the second part of Theorem \ref{pml} is concerned, basically speaking, through very involved, it follows from the first part and the fact that the approximants $\xi_{k,n}$ of the continued fraction related to $A_k\cdots A_n$ are convergent in an exponential rates, see Lemma \ref{dxr} below.

  We remark that whenever $\textrm{det}(M_k)>0,$ $A_k$'s  constructed  in \eqref{dta} are not nonnegative matrices any longer, so some special techniques are required when considering the entries of their products.

{\bf Outline of the paper.} The remainder of the paper is arranged as follows. In Section \ref{pr}, we give some preliminary results of continued fractions and list some basic facts of the matrices $M_k$ and $A_k,k\ge1.$ In Section \ref{pc},  we estimate firstly the entries of $A_{k+1}\cdots A_{k+n}$ and $\sum_{s=1}^{n+1}A_{k+s}\cdots A_{k+n}$ and then transit to the asymptotics of the entries of $M_{k+1}\cdots M_{k+n}$ and $\sum_{s=1}^{n+1}M_{k+s}\cdots M_{k+n}.$

\section{Preliminary results}\label{pr}
 Products of 2-by-2 matrices are closely related to continued fractions and therefore the theory of continued fractions plays a key role in proving Theorem \ref{pml}. To begin with, we introduce some basics of continued fractions.
 \subsection{Continued fractions and their tails}

Let $\beta_k,\alpha_k,k\ge 1$ be certain  real numbers. For $1\le k\le n,$ We denote by
\begin{equation}\label{aprx}
\xi_{k,n}\equiv\frac{\beta_k}{\alpha_k}\begin{array}{c}
                                \\
                               +
                             \end{array}\frac{\beta_{k+1}}{ \alpha_{k+1}}\begin{array}{c}
                                \\
                               +\cdots+
                             \end{array}\frac{\beta_n}{\alpha_n}:=\dfrac{\beta_k}{\alpha_k+\dfrac{\beta_{k+1}}{\alpha_{k+1}+_{\ddots_{\textstyle +\frac{\textstyle\beta_{n}}{\textstyle\alpha_{n} } }}}}
\end{equation}
 the $(n-k+1)$-th approximant of a continued fraction
 \begin{align}\label{xic}
   &\xi_k:=\frac{\beta_{k}}{\alpha_{k}}\begin{array}{c}
                                \\
                               +
                             \end{array}\frac{\beta_{k+1}}{\alpha_{k+1 }}\begin{array}{c}
                                \\
                               +
                             \end{array}\frac{\beta_{k+2}}{\alpha_{k+2}}\begin{array}{c}
                                \\
                               +\cdots
                             \end{array}.
\end{align}
In the literature, $\xi_k,k\ge1$ in \eqref{xic}  are usually called the tails
of the continued fraction ${\frac{\beta_{1}}{\alpha_{1}}}_{+}\frac{\beta_{2}}{\alpha_{2 }}_{+\cdots}.$
 If  $\lim_{n\rightarrow\infty}\xi_{k,n}$ exists, then we say that the continued fraction $\xi_k$  is convergent and its value is defined as $\lim_{n\rightarrow\infty}\xi_{k,n}.$  If $\beta_k>0,\alpha_k>0, \forall k\ge1$ and \begin{align}
  \exists C>0 \text{ such that } \forall k\ge1,\ C^{-1}\le {\beta_k}/{\alpha_k}\le C,\label{ssc}
 \end{align}
 then  by Seidel-Stern Theorem (see Lorentzen and Waadeland \cite[Theorem 3.14]{lw}), for any $k\ge1,$  $\xi_k$ is convergent.
The lemma below will also be  used times and again.
 \begin{lemma}\label{ct}
If $
  \lim_{n\rto}\alpha_n=\alpha\ne0,$ $\lim_{n\rto}\beta_n=\beta$  and  $\alpha^2+4\beta\ge0,$
 then
 for any $k\ge1,$  $\lim_{n\rto}\xi_{k,n}$ exists and furthermore \begin{align}
  \lim_{k\rto}\xi_k=\frac{\alpha}{2}\z(\sqrt{1+4\beta/\alpha^2}-1\y).\no
 \end{align}
 \end{lemma}
The proof Lemma \ref{ct} can be found in many references,
  we refer the reader to \cite{lor}, see  discussion between (4.1) and (4.2) on page 81 therein.

The following lemma gives various inequalities related to the tails and theirs approximants of a continued fraction which are useful for us.
 \begin{lemma}\label{xp} Let $\xi_{k,n}$ and $\xi_k$ be the ones in \eqref{aprx} and \eqref{xic}. Suppose that $\alpha_k,\beta_k>0, \forall k\ge1$ and \eqref{ssc} is satisfied.
 Then we have $\xi_{k,n}\rightarrow\xi_k\in(0,\infty), \text{ as } n\rightarrow\infty$ and
 \begin{align} 
  & \xi_{k,n}\Big\{\begin{array}{cc}
     <\xi_k,& \text{if } n-k+1 \text{ is even,} \\
     >\xi_k,&\text{if } n-k+1 \text{ is odd,}
   \end{array} 1\le k\le n,\label{pxi}\\
   &\xi_{k,n}\xi_{k+1,n}\Big\{\begin{array}{cc}
     >\xi_k\xi_{k+1},& \text{if } n-k+1 \text{ is even,} \\
     <\xi_k\xi_{k+1},&\text{if } n-k+1 \text{ is odd,}
   \end{array}1\le k\le n-1. \label{txi}
 \end{align}
    \end{lemma}
    For the proof of the lemma, we refer the reader to \cite[Lemma 5.3]{wy}.
\subsection{Some facts on  matrices $A_k$ and $M_k$}\label{ma}

Let $M_k,k\ge1$ be those in \eqref{mg} and $A_k,\Lambda_k,k\ge1$ be those in \eqref{am}.  In this section we always assume that $$\lim_{k\rto}M_k=M,   a+\theta\ne0, b\ne0\text{ and }bd\ne a\theta.$$
Then it is easily seen that \begin{equation}\no
 \lim_{k\rto} A_k=A:=\z(\begin{array}{cc}
                          a+\theta & b \\
                          d-a\theta/b & 0
                        \end{array}
  \y)
\end{equation}
and clearly $A$ is similar to $M$ so that  they share the same eigenvalues. Furthermore, we can find $\varepsilon>0$ and $k_1>0$ such that
\begin{align}
   &\label{atg}a_k+\theta_k>\varepsilon, b_k>\varepsilon, \forall k\ge k_1; \\
    &\label{bda}bd>a\theta\Rightarrow   \tilde d_k>\varepsilon, \forall k\ge k_1;\\
     &\label{bdb}bd<a\theta\Rightarrow  \tilde d_k<\varepsilon, \forall k\ge k_1.
\end{align}
Therefore, using \eqref{atg}, since $A_k\cdots A_n=\Lambda_k^{-1}M_k\cdots M_n\Lambda_{n+1},n\ge k\ge1,$  we must have
\begin{align}
    \mb e_1 A_k\cdots A_n \mb e_1^t&=\mb e_1M_k\cdots M_n(1,\theta_{n+1}/b_{n+1})^t>0,n>k\ge k_1,\label{apm1}\\
  \mb e_1 A_k\cdots A_n \mb e_2^t&=\mb e_1M_k\cdots M_n\mb e_2^t>0,n> k\ge k_1.\label{apm2}
    \end{align}
Finally, taking \eqref{bda}-\eqref{apm2} into account, we get
 \begin{equation*}\begin{split}
   &bd<a\theta\Rightarrow \mb e_2 A_k\cdots A_n \mb e_1^t<0,\mb e_2 A_k\cdots A_n \mb e_2^t<0, n\ge k\ge k_1;\\
    &bd>a\theta \Rightarrow \mb e_2 A_k\cdots A_n \mb e_1^t>0,\mb e_2 A_k\cdots A_n \mb e_2^t>0, n\ge k\ge k_1.\\
    \end{split}
\end{equation*}

\section{Product of matrices and continued fractions}\label{pc}
As explained in Section \ref{intr}, it is hard to study directly the asymptotics of entries of $M_{k+1}\cdots M_{k+n}$ as $n\rto.$ In this section, we will treat the entries of $A_{k+1}\cdots A_{k+n}$ by using the theory of continued fractions.  The main task of this section is to prove the  following theorem which leads to Theorem \ref{pml}.
\begin{theorem}\label{s}
 Suppose that $\lim_{k\rto}M_k=M,$   $a+\theta\ne0,$ $b\ne0$ and $bd\ne a\theta.$ Then $\exists k_0>0$ such that for $k\ge k_0$ and $i,j\in\{1,2\},$  we have
  \begin{align}\label{pp}
    &\quad\ \ \frac{\mb e_iA_{k+1}\cdots A_{k+n}\mb e_j^t}{\xi_{k+1}^{-1}\cdots\xi_{k+n}^{-1}}\rightarrow\phi(i,j,k),
    \end{align} uniformly in $k$ as $n\rto,$ where $\phi(1,1,k)=\frac{\varrho}{\varrho-\varrho_1},$
  $\phi(1,2,k)=\frac{b}{\varrho-\varrho_1},$
  $\phi(2,1,k)=\frac{\tilde d_{k+1}}{\xi_{k+1}^{-1}}\frac{\varrho}{\varrho-\varrho_1}$
and  $\phi(2,2,k)=\frac{\tilde d_{k+1}}{\xi_{k+1}^{-1}}\frac{b}{\varrho-\varrho_1}.$ Moreover, if we assume further $\varrho\ge1,$ then for $k\ge k_0$ and $i=1,$ $j\in \{1,2\},$ with the above $\phi(i,j,k),$ we have
    \begin{align}
    \label{ps}
    &\frac{\sum_{s=1}^{n+1}\mb e_iA_{k+s}\cdots A_{k+n}\mb e_j^t}{\sum_{s=1}^{n+1}\xi_{k+s}^{-1}\cdots \xi_{k+n}^{-1}}\rightarrow\phi(i,j,k),
  \end{align}
   uniformly in $k$ as $n\rto.$
\end{theorem}

\subsection{Product of matrices expressed in terms of the approximants of continued fractions}
In what follows,  $k_1$ will be always the one in Section \ref{ma}. We first write $\mb e_1A_k\cdots A_n\mb e_1^t$ in terms of the approximants of some continued fractions.   For $k_1\le k\le n,$ we set
 \begin{align}\label{xiy}
   y_{k,n}:=\mathbf e_1 A_k\cdots A_{n}\mathbf e_1^t \text{ and }\xi_{k,n}:=\frac{y_{k+1,n}}{y_{k,n}}.
 \end{align}
Then using the fact $y_{n+1,n}=1,$ we get
 \begin{align}
\xi_{k,n}^{-1}\cdots \xi_{n,n}^{-1}&=y_{k,n}=\mathbf e_1 A_k\cdots A_{n}\mathbf e_1^t.\label{cpn}
 \end{align}
 \begin{lemma}\label{axc}For $k_1\le k\le n,$ $\xi_{k,n}$ defined in \eqref{xiy} coincides with the one in \eqref{aprx} with $\beta_k=\tilde b_k^{-1}\tilde d_{k+1}^{-1}$ and $\alpha_k=\tilde a_k\tilde b_k^{-1}\tilde d_{k+1}^{-1}.$
    \end{lemma}
    \proof  Clearly, $\xi_{n,n}=\frac{1}{y_{n,n}}=\frac{1}{\tilde a_n}=\frac{\tilde b_n^{-1}\tilde d_{n+1}^{-1}}{\tilde a_n\tilde b_n^{-1}\tilde d_{n+1}^{-1}}.$ For $k_1\le k< n,$ note that
\begin{align}\label{ix}
  \xi_{k,n}&=\frac{y_{k+1,n}}{y_{k,n}}=\frac{\mathbf e_1 A_{k+1}\cdots A_{n}\mathbf e_1^t}{\mathbf e_1 A_k\cdots A_{n}\mathbf e_1^t}=\frac{\mathbf e_1 A_{k+1}\cdots A_{n}\mathbf e_1^t}{(\tilde a_k\mathbf e_1+\tilde b_k\mathbf e_2) A_{k+1}\cdots A_{n}\mathbf e_1^t}\\
  &=\frac{1}{\tilde a_k+\tilde b_k\frac{\mathbf e_2 A_{k+1}\cdots A_{n}\mathbf e_1^t}{\mathbf e_1 A_{k+1}\cdots A_{n}\mathbf e_1^t}}=\frac{1}{\tilde a_k+\tilde b_k\tilde d_{k+1}\frac{\mathbf e_1 A_{k+2}\cdots A_{n}\mathbf e_1^t}{\mathbf e_1 A_{k+1}\cdots A_{n}\mathbf e_1^t}}\no\\
  &=\frac{\tilde b_k^{-1}\tilde d_{k+1}^{-1}}{\tilde a_k\tilde b_k^{-1}\tilde d_{k+1}^{-1}+\xi_{k+1,n}}.\no
  \end{align}
  Thus, the lemma can be proved by iterating \eqref{ix}. \qed
\subsection{Proof of Theorem \ref{s}}

 We give the proof of Theorem \ref{s} in this section. To begin with, we derive some auxiliary lemmas.  Keep in mind that in what follows, unless otherwise specified, we always assume that $c>0$ is a constant which may change from line to line.
\subsubsection{Auxiliary lemmas}
 Let $\xi_k,\xi_{k,n}, n\ge k_1$ be those defined in \eqref{aprx} and \eqref{xic} with $\beta_k=\tilde b_k^{-1}\tilde d_{k+1}^{-1}$ and $\alpha_k=\tilde a_k\tilde b_k^{-1}\tilde d_{k+1}^{-1}.$  Suppose that the conditions of Theorem \ref{s} are all fulfilled. Then we have \begin{align}\label{bal}
  \lim_{k\rto}\beta_k=:\beta=(bd-a\theta)^{-1}\ne 0\text{ and }\lim_{k\rto}\alpha_k=:\alpha=\frac{a+\theta}{bd-a\theta}\ne 0.
\end{align}
Noticing that $\alpha^2+4\beta=\frac{(a-\theta)^2+4bd}{(bd-a\theta)^2}\ge0,$ it thus
 follows from Lemma \ref{ct} that
\begin{align}\label{lxk}
  \lim_{n\rto}\xi_{k,n}=\xi_k\ \text{\ exists and } \lim_{k\rto}\xi_k=:\xi=\frac{\alpha}{2}\z(\sqrt{1+4\beta/\alpha^2}-1\y)=\varrho^{-1}>0.
\end{align}
Moreover, by  \eqref{apm1} and \eqref{xiy} we have
\begin{align}\label{po}
 \xi_k>0, \xi_{k,n}>0,\ \forall n\ge k\ge k_1.
\end{align}
Therefore, from \eqref{lxk} and \eqref{po} we get
\begin{align}\label{up}
  c^{-1}<\xi_k<c \text{ and } c^{-1}<\xi_{k,n}<c, n\ge k\ge k_1.
\end{align}
We remark that $\xi_n, n\ge k_1$ here coincide with those defined in \eqref{xicn}. The following two auxiliary lemmas play important roles in the proof of Theorem \ref{s}.
\begin{lemma}\label{mi}Suppose that  all conditions of Theorem \ref{s} are fulfilled and $bd<a\theta.$  Fix $k\ge k_1.$ Then $\xi_{k,n},n\ge k$ is monotone increasing in $n$ and thus $\xi_{k,n}<\xi_k,n\ge k.$
\end{lemma}

\begin{lemma}\label{dxr} Suppose that  all conditions of Theorem \ref{s} are fulfilled.  Then $\exists 0<r<1$ and $ k_2>0$ such that
\begin{align}\label{xd}
 |\xi_{k}-\xi_{k,n}|\le r^{n-k}|\xi_n-\xi_{n,n}|, \forall n\ge k>k_2.
\end{align}
\end{lemma}

\noindent{\it Proof of Lemma \ref{mi}. } Write temporarily $A_{k,n}:=A_k\cdots A_n$ for $n\ge k\ge1.$ Some easy computation yields that
\begin{align}\label{x}
  \xi_{k,n+1}&-\xi_{k,n}=\frac{A_{k+1,n+1}(11)}{A_{k,n+1}(11)}-\frac{A_{k+1,n}(11)}{A_{k,n}(11)}\\
  &=\frac{A_{k,n}(11)A_{k+1,n+1}(11)-A_{k+1,n}(11)A_{k,n+1}(11)}{A_{k,n}(11)A_{k,n+1}(11)}\no\\
  &=-\tilde b_k\tilde d_{n+1}\frac{A_{k+1,n}(11)A_{k+1,n}(22)-A_{k+1,n}(12)A_{k+1,n}(21)}{A_{k,n}(11)A_{k,n+1}(11)}\no\\
  &=\frac{-\tilde b_k\tilde d_{n+1}\text{det}(A_{k+1,n})}{A_{k,n}(11)A_{k,n+1}(11)}=\frac{-\tilde b_k\tilde d_{n+1}\prod_{j=k+1}^n\text{det}(A_j)}{A_{k,n}(11)A_{k,n+1}(11)}.\no
    \end{align}
  By assumption, for $k\ge k_1$ we must have  $\tilde a_k>0,$ $\tilde b_k>0,$ $\tilde d_{k}<0$ and therefore $\mathrm{det}(A_k)>0$ and $A_{k,n}(11)>0,n\ge k\ge k_1.$ Consequently, we conclude from \eqref{x} that  $\xi_{k,n+1}-\xi_{k,n}>0,n\ge k\ge k_1.$ The lemma is proved. \qed

\noindent{\it Proof of Lemma \ref{dxr}.}
Assume first $bd>a \theta.$  We claim that $\exists k_2>k_1$ such that
\begin{align}
  n-k+1\text{ is even } &\Rightarrow \frac{\xi_{k,n}-\xi_k}{\xi_{k+1}-\xi_{k+1,n}}< r, \forall n>k\ge k_2, \label{udo}\\
  n-k+1\text{ is odd } &\Rightarrow \frac{\xi_{k,n}-\xi_k}{\xi_{k+2,n}-\xi_{k+2}}<r^2,\forall n>k\ge k_2,\label{udt}
\end{align}
where $0<r<1$ is a proper number.

 In fact, by \eqref{bda}, we have $\beta_k>0,\alpha_k>0, \forall k\ge k_1$ and
  \begin{equation}\label{ab}
    \alpha_k\rightarrow \alpha>0, \beta_k\rightarrow\beta>0, \text{ as }k\rto.
  \end{equation}
 It then follows from \eqref{lxk} that
 $\frac{\xi_{k}}{\alpha_k+\xi_{k+1}}\rightarrow\frac{\xi}{\alpha+\xi}<1$ as $k\rightarrow\infty.$
 As a result,  for some proper number $0<r<1,$ $\exists k_2>k_1$ such that
 $\frac{\xi_{k}}{\alpha_k+\xi_{k+1}}<r, \forall k\ge k_2.$
 On the other hand,
 it follows from \eqref{aprx} and \eqref{xic} that
 \begin{align*}
 \frac{\xi_{k,n}-\xi_k}{\xi_{k+1}-\xi_{k+1,n}}=\frac{\xi_{k,n}}{\alpha_k+\xi_{k+1}}=\frac{\xi_{k}}{\alpha_k+\xi_{k+1,n}}.
 \end{align*}
 In view of \eqref{ab}, we can also apply Lemma 2 to $\xi_{k,n},\xi_k,k\ge k_1.$
 If $n-k+1$ is even, then by \eqref{pxi}, $\xi_{k+1,n}>\xi_{k+1},n> k\ge k_1.$ Thus
 \begin{align*}
 \frac{\xi_{k,n}-\xi_k}{\xi_{k+1}-\xi_{k+1,n}}=\frac{\xi_{k}}{\alpha_k+\xi_{k+1,n}}<\frac{\xi_{k}}{\alpha_k+\xi_{k+1}}<r, \ \forall n>k>k_2.
 \end{align*}
  If  $n-k+1$ is odd, then by \eqref{txi}, we have $\xi_{k,n}\xi_{k+1,n}<\xi_k\xi_{k+1},n> k\ge k_1.$ Therefore,
 \begin{align*}
 \frac{\xi_{k,n}-\xi_k}{\xi_{k+2,n}-\xi_{k+2}}
 &=\frac{\xi_{k,n}}{\alpha_k+\xi_{k+1}}\frac{\xi_{k+1,n}}{\alpha_{k+1}+\xi_{k+2}}\\
 &<\frac{\xi_{k}}{\alpha_k+\xi_{k+1}}\frac{\xi_{k+1}}{\alpha_{k+1}+\xi_{k+2}}<r^2,\ \forall n\ge k\ge k_2.
 \end{align*}
 We thus finish  proving the claim and \eqref{xd} is a direct consequence of \eqref{udo} and \eqref{udt}.

Next we  assume $bd<a\theta.$ Then from \eqref{bdb}, we get $\tilde d_k\le -\varepsilon,\ \forall k\ge k_1.$
 Taking Lemma \ref{mi} into account, since $\xi_{k,n}>0,$ $\forall n\ge k_1,$ then by some easy computation, we have
 from \eqref{aprx} and \eqref{xic} that
 \begin{align*}
 0<\frac{\xi_{k,n}-\xi_k}{\xi_{k+1,n}-\xi_{k+1}}=-\frac{\xi_{k,n}}{\alpha_k+\xi_{k+1}}<-\frac{\xi_{k}}{\alpha_k+\xi_{k+1}}.
 \end{align*}
But
 $-\frac{\xi_{k}}{\alpha_k+\xi_{k+1}}\rightarrow-\frac{\xi}{\alpha+\xi}=\frac{\alpha+\sqrt{\alpha^2+4\beta}}{\alpha-\sqrt{\alpha^2+4\beta}}<1$ as $k\rto.$
 As a result,  for some proper number $0<r<1,$ $\exists k_2>0$ such that
 $-\frac{\xi_{k}}{\alpha_k+\xi_{k+1}}<r, \forall k> k_2,$ which finishes the proof of the lemma. \qed

%
%

\subsubsection{Proof of \eqref{pp}}
To begin with, we deal with the case $i=j=1.$ In what follows we set $k_0=k_1\vee k_2.$ Write $\Pi(k,n)=\frac{\xi_{k+1}^{-1}\cdots\xi_{k+n}^{-1}}{\mb e_1A_{k+1}\cdots A_{k+n}\mb e_1^t},n\ge 1,k\ge k_0.$ Then it follows from \eqref{cpn}
that
\begin{align*}
  \Pi(k,n)=\frac{\xi_{k+1}^{-1}\cdots\xi_{k+n}^{-1}}{\xi_{k+1,k+n}^{-1}\cdots \xi_{k+n,k+n}^{-1}}=\frac{\xi_{k+1,k+n}\cdots \xi_{k+n,k+n}}{\xi_{k+1}\cdots\xi_{k+n}}
\end{align*}
 Let $\overline{\pi}(k):=\varlimsup_{n\rto}\log \Pi(k,n)$ and $\underline{\pi}(k)=\varliminf_{n\rto}\log \Pi(k,n).$
Since $\lim_{n\rto}\xi_n=\varrho^{-1}$ and for each $i>0,$ $\lim_{j\rto}A_{n+1} A_{n+2}\cdots A_{n+i}=A^i,$ we have
\begin{align*}
  \lim_{n\rto}&\log\frac{\xi_{k+n-i+1,k+n}\cdots \xi_{k+n,k+n}}{\xi_{k+n-i+1}\cdots\xi_{k+n}}=\lim_{n\rto}\log\frac{\z(\mathbf e_1A_{k+n-i+1}\cdots A_{k+n}\mathbf e_1^t\y)^{-1}}{\xi_{k+n-i+1}\cdots\xi_{k+n}}\\
  &=\log\z(\varrho^{-i}\mathbf e_1A^i\mathbf e_1^t\y)^{-1}=\log\frac{\varrho^i(\varrho-\varrho_1)}{\varrho^{i+1}-\varrho_1^{i+1}}=:\sigma_i,
\end{align*} where  the convergence is uniform in $k.$
Then, consulting to \eqref{up} and \eqref{xd}, we have
\begin{align*}
  \overline{\pi}(k)&=\varlimsup_{n\rto}\log \Pi(k,n)\\
  &=\varlimsup_{n\rto}\z(\sum_{j=1}^{n-i}\log\z(1+\frac{\xi_{k+j,k+n}-\xi_{k+j}}{\xi_{k+j}}\y)+
  \log\frac{\xi_{k+n-i+1,k+n}\cdots \xi_{k+n,k+n}}{\xi_{k+n-i+1}\cdots\xi_{k+n}}\y)\\
  &=\sigma_i+\varlimsup_{j\rto}\sum_{j=1}^{n-i}\log\z(1+\frac{\xi_{k+j,k+n}-\xi_{k+j}}{\xi_{k+j}}\y)\\
  &\le \sigma_i+c\varlimsup_{n\rto}\sum_{j=1}^{n-i}|\xi_{k+j,k+n}-\xi_{k+j}|\\
  &\le \sigma_i+c\varlimsup_{n\rto}\sum_{j=1}^{n-i}r^{n-j}=\sigma_i+c\varlimsup_{n\rto}\sum_{j=i}^{n-1}r^{j}= \sigma_i+cr^i.
\end{align*}

For a lower limit, note that by \eqref{xd}, we have
\begin{align}\label{ci}
  |\xi_{k+j,k+n}-\xi_{k+j}|/\xi_{k+j}<r^{n-j}|\xi_{k+n,k+n}-\xi_{k+n}|/\xi_{k+j}<cr^{n-j}.
\end{align}
Fix $i_0$ such that $cr^{i_0}<1,$ where $c$ is the one in \eqref{ci}.
Then for $i\ge i_0,$
using the fact $\log(1-x)\ge -x/(1-x),1>x\ge 0,$ we have
\begin{align*}
  \underline{\pi}(k)&=\varliminf_{n\rto}\log \Pi(k,n)\\
  &=\varliminf_{n\rto}\z(\sum_{j=1}^{n-i}\log\z(1+\frac{\xi_{k+j,k+n}-\xi_{k+j}}{\xi_{k+j}}\y)+
  \log\frac{\xi_{k+n-i+1,k+n}\cdots \xi_{k+n,k+n}}{\xi_{k+n-i+1}\cdots\xi_{k+n}}\y)\\
  &=\sigma_i+\varliminf_{j\rto}\sum_{j=1}^{n-i}\log\z(1+\frac{\xi_{k+j,k+n}-\xi_{k+j}}{\xi_{k+j}}\y)\\
  &\ge \sigma_i-c\varliminf_{n\rto}\sum_{j=1}^{n-i}|\xi_{k+j,k+n}-\xi_{k+j}|\\
  &\ge \sigma_i-c\varliminf_{n\rto}\sum_{j=1}^{n-i}r^{n-j}=\sigma_i-c\varliminf_{n\rto}\sum_{j=i}^{n-1}r^{j}= \sigma_i-cr^i.
\end{align*}
We have thus shown that for each $i\ge i_0,$
\begin{align}\label{plu}
  \sigma_i-cr^i=\varliminf_{n\rto}\log \Pi(k,n)\le \varlimsup_{n\rto}\log \Pi(k,n)\le\sigma_i+cr^i
\end{align} where the upper and lower limits we take here are indeed uniform in $k.$
Now we can infer that
\begin{align*}
  \overline\pi(k)-\underline\pi(k)\le cr^i,i\ge i_0.
\end{align*}
Letting $i\rto,$ we get $\overline\pi(k)=\underline\pi(k)$
so that for each $k\ge k_0,$ the limit
$$\pi(k):=\lim_{n\rto}\log \Pi(k,n)$$ exists  and moreover the convergence is uniform in $k.$
 Since $\lim_{n\rto}\xi_{k+j,k+n}=\xi_{k+j},k\ge k_0,j\ge 1$ we conclude that $\pi(k)$ is independent of $k.$
 In addition, from \eqref{plu} we have
 $|\pi(k)-\sigma_i|\le cr^i,k\ge k_0, i\ge i_0.$  Letting $i\rto$ we get
 \begin{align*}
   \pi(k)\equiv \lim_{i\rto}\sigma_i=\log\frac{\varrho-\varrho_1}{\varrho}, k\ge k_0.
 \end{align*}
 We can now come to the conclusion that
 $$\lim_{n\rto}\Pi(k,n)=\frac{\varrho-\varrho_1}{\varrho}, k\ge k_0$$ and the convergence is uniform in $k,$ which finishes the proof of the case $i=j=1.$

 Next, we consider the case $i=j=2.$ Notice that
 \begin{align*}
   \frac{\mb e_2A_{k+1}\cdots A_{k+n}\mb e_2^t}{\xi_{k+1}^{-1}\cdots\xi_{k+n}^{-1}}=\frac{\tilde d_{k+1}}{\xi_{k+1}^{-1}} \frac{\mb e_1A_{k+2}\cdots A_{k+n-1}\mb e_1^t}{\xi_{k+2}^{-1}\cdots\xi_{k+n-1}^{-1}}\frac{\tilde b_{k+n}}{\xi_{k+n}^{-1}}.
 \end{align*}
 Since $\lim_{n\rto} \tilde b_{k+n}=b$ and $\lim_{n\rto}\xi_{n+k}=\varrho^{-1},$ it then follows from the case $i=j=1$ that uniformly in $k,$
 $$ \lim_{n\rto}\frac{\mb e_2A_{k+1}\cdots A_{k+n}\mb e_2^t}{\xi_{k+1}^{-1}\cdots\xi_{k+n}^{-1}}=\frac{\tilde d_{k+1}}{\xi_{k+1}^{-1}}\frac{b}{\varrho-\varrho_1},k\ge k_0.$$
 Since the proof of the case $i=1,j=2$ and that of the case $i=2,j=1$ can be given by some similar arguments,  we can complete the proof here.
  \qed

  \subsubsection{Proof of \eqref{ps}}
  We prove first the case $i=j=1.$ Taking \eqref{cpn} into account, we obtain
  \begin{align*}
    &\frac{\sum_{s=1}^{n+1}\mb e_1A_{k+s}\cdots A_{k+n}\mb e_1^t}{\sum_{s=1}^{n+1}\xi_{k+s}^{-1}\cdots \xi_{k+n}^{-1}}=\frac{\sum_{s=1}^{n+1} \xi_{k+s,k+n}^{-1}\cdots \xi_{k+n,k+n}^{-1}}{\sum_{s=1}^{n+1}\xi_{k+s}^{-1}\cdots \xi_{k+n}^{-1}}\\
    &\quad\quad\quad\quad\quad\quad=\frac{\sum_{s=1}^{n+1} \xi_{k+1,k+n}\cdots \xi_{k+s-1,k+n}}{\sum_{s=1}^{n+1}\xi_{k+1}\cdots \xi_{k+s-1}}\frac{\xi_{k+1}\cdots\xi_{k+n}}{\xi_{k+1,k+n}\cdots\xi_{k+n,k+n}}\\
    &\quad\quad\quad\quad\quad\quad =\frac{\sum_{s=1}^{n+1} \xi_{k+1,k+n}\cdots \xi_{k+s-1,k+n}}{\sum_{s=1}^{n+1}\xi_{k+1}\cdots \xi_{k+s-1}}\frac{\mb e_1A_{k+1}\cdots A_{k+n}\mb e_1^t}{\xi_{k+1}^{-1}\cdots\xi_{k+n}^{-1}}.
  \end{align*}
  In view of \eqref{pp}, in order to prove \eqref{ps}, it suffices to show that  for $k\ge k_0,$
  \begin{align}\label{cone}
    \frac{\sum_{s=1}^{n+1} \xi_{k+1,k+n}\cdots \xi_{k+s-1,k+n}}{\sum_{s=1}^{n+1}\xi_{k+1}\cdots \xi_{k+s-1}}\rightarrow 1,
  \end{align} uniformly in $k$ as $n\rto.$

Assume first $bd<a\theta.$ It then follows from Lemma \ref{mi} that $\xi_{k,n}<\xi_k$ for $n\ge k\ge k_0$ and consequently we get
\begin{align}\label{nu}
  \sum_{s=1}^{n+1} \xi_{k+1,k+n}\cdots \xi_{k+s-1,k+n}< \sum_{s=1}^{n+1}\xi_{k+1}\cdots\xi_{k+s-1},n\ge1, k\ge k_0.
\end{align}
For a lower bound, applying the lemmas \ref{mi} and \ref{dxr} and consulting to \eqref{up}, we obtain
\begin{align}\label{it}
  &\sum_{s=1}^{n+1} \xi_{k+1,k+n}\cdots \xi_{k+s-1,k+n}=1+\sum_{s=1}^{n} \xi_{k+1,k+n}\cdots \xi_{k+s,k+n}\\
  &\quad\quad=1+\xi_{k+1}+\xi_{k+1}\sum_{s=2}^{n} \xi_{k+2,k+n}\cdots \xi_{k+s,k+n}\no\\
  &\quad\quad\quad\quad+(\xi_{k+1,k+n}-\xi_{k+1})\z(1+\sum_{s=2}^{n} \xi_{k+2,k+n}\cdots \xi_{k+s,k+n}\y)\no\\
  &\quad\quad\ge 1+\xi_{k+1}+\xi_{k+1}\sum_{s=2}^{n} \xi_{k+2,k+n}\cdots \xi_{k+s,k+n}\no\\
  &\quad\quad\quad\quad-cr^{n-k}\z(1+\sum_{s=2}^{n} \xi_{k+2}\cdots \xi_{k+s}\y),k\ge k_0,n\ge1\no.
\end{align}
Iterating \eqref{it},  we have
\begin{align}\label{la}
\sum_{s=1}^{n+1} \xi_{k+1,k+n}&\cdots \xi_{k+s-1,k+n}\ge \sum_{s=1}^{n+1} \xi_{k+1}\cdots \xi_{k+s-1}\\
  &-c\sum_{i=1}^nr^{n-i}\xi_{k+1}\cdots\xi_{k+i-1}\z(1+\sum_{s=i+1}^{n}\xi_{k+i+1}\cdots\xi_{k+s}\y)\no\\
  &=:\mathrm{(I)}-\mathrm{(II)}, k\ge k_0,n\ge1.\no
\end{align}
Clearly, for $1\le i\le n,$
\begin{align}\label{lo}
  \frac{\xi_{k+1}\cdots\xi_{k+i}\z(1+\sum_{s=i+1}^{n}\xi_{k+i+1}\cdots\xi_{k+s}\y)}{\sum_{s=1}^{n+1}\xi_{k+1}\cdots\xi_{k+s-1}}\le 1,k\ge k_0,n\ge1.
\end{align}
Let $\varepsilon>0$ be an arbitrary number and fix $k_3$ such that $r^{k_3}<\varepsilon.$
We claim that for $k\ge k_0,$
 \begin{align}\label{uc}
   \frac{\xi_{k+1}\cdots\xi_{k+n-i}}{\sum_{s=1}^{n+1} \xi_{k+1}\cdots \xi_{k+s-1}}\rightarrow 0
 \end{align} uniformly in $k$ and $0\le i\le k_3$ as $n\rto.$

 To prove the claim,
 recall that by \eqref{lxk} and \eqref{po}, $\xi_k>0,k\ge k_0$ and $\lim_{k\rto}\xi_k^{-1}=\varrho\ge1.$ Suppose first $\varrho>1.$ Then for certain $1<\sigma<\varrho,$ there exists $N_0>0$ such that $\xi_n<\sigma^{-1}, \forall n\ge N_0.$  Therefore using \eqref{up},  for $k\ge k_0,$ we obtain that
\begin{align*}
\frac{\xi_{k+1}\cdots\xi_{k+n-i}}{\sum_{s=1}^{n+1} \xi_{k+1}\cdots \xi_{k+s-1}}\le c^{N_0}\xi_{k+N_0+1}\cdots\xi_{k+n-i}\le c^{N_0}\sigma^{-(n-N_0-i)}\rightarrow 0,
\end{align*}
 uniformly in $k$ and $0\le i\le k_3$ as $n\rto.$ Suppose next $\varrho=1.$ Fix an integer $M\ge1$ and a number $1>\sigma_1>0.$ There exists a number $N_1>1$ independent of $k$ such that for $n>N_1,$ we have $\xi_{k+n-i+1}^{-1}\cdots\xi_{k+n}^{-1}<1+\sigma_1, 0\le i\le k_3$ and $\xi_{k+n-j}^{-1}\cdots \xi_{n+k}^{-1}>1-\sigma_1, \forall 0\le j\le M.$
 Consequently, for $0\le i\le k_3, k\ge k_0,$ we have
 \begin{align*}
\frac{\xi_{k+1}\cdots\xi_{k+n-i}}{\sum_{s=1}^{n+1} \xi_{k+1}\cdots \xi_{k+s-1}}&=\frac{\xi_{k+n-i+1}^{-1}\cdots\xi_{k+n}^{-1}}{\sum_{s=1}^{n+1}\xi_{k+s}^{-1}\cdots\xi_{k+n}^{-1}}
   \le \frac{\xi_{k+n-i+1}^{-1}\cdots\xi_{k+n}^{-1}}{\sum_{s=n-M}^{n+1}\xi_{k+s}^{-1}\cdots\xi_{k+n}^{-1}}\\
  & <(1+\sigma_1)(1-\sigma_1)^{-1}(M+2)^{-1}.
 \end{align*}
  Since $M$ is arbitrary, we complete the proof of the claim.

 Now, we are ready to deal with the term (II) on the right-hand side of \eqref{la}. Taking \eqref{up}, \eqref{lo} and \eqref{uc} into consideration, for $k\ge k_0,$ we have
 \begin{align*}
   \varlimsup_{n\rto}&\frac{\mathrm{(II)}}{\sum_{s=1}^{n+1} \xi_{k+1}\cdots \xi_{k+s-1}}\le c\varlimsup_{n\rto} \sum_{i=1}^{n-k_3}r^{n-i}\\
   &   +c \varlimsup_{n\rto} \sum_{i=n-k_3+1}^{n}r^{n-i}\frac{\xi_{k+1}\cdots\xi_{k+i}\z(1+\sum_{s=i+1}^{n}\xi_{k+i+1}\cdots\xi_{k+s}\y)}{\sum_{s=1}^{n+1}\xi_{k+1}\cdots\xi_{k+s-1}}\\
   &=c\varlimsup_{n\rto} \sum_{i=1}^{n-k_3}r^{n-i}=cr^{k_3}/(1-r)\le c\varepsilon/(1-r),
 \end{align*} where we emphasize that the second term in the middle converges to $0$ uniformly in $k.$  Since $\varepsilon$ is arbitrary, for $k\ge k_0$ we have
 \begin{align}\label{lb}
   \frac{\mathrm{(II)}}{\sum_{s=1}^{n+1} \xi_{k+1}\cdots \xi_{k+s-1}}\rightarrow 0
 \end{align} uniformly in $k$ as $n\rto.$ Taking \eqref{nu}, \eqref{la} and \eqref{lb} together, we get \eqref{cone} whenever $bd<a\theta.$

Next, we assume $bd>a\theta.$  Then taking \eqref{bda} into account, we have $\beta_k=\tilde b_k^{-1}\tilde d_{k+1}^{-1}>0$ and $\alpha_k=\tilde a_k\tilde b_k^{-1}\tilde d_{k+1}^{-1}>0,$  for $k\ge k_0$ and by \eqref{bal} we get $c^{-1}\le \max\{\beta_k,\alpha_k\}\le c, k\ge k_0.$
 Thus the conditions of  Lemma \ref{xp} are fulfilled. We consider below only the case $n$ is even since the case $n$ is odd follows similarly. Applying Lemma \ref{xp}, we obtain
  \begin{align}\label{fu}
 \sum_{s=1}^{n+1}\prod_{j=k+1}^{k+s-1}\xi_{j,k+n}\le \sum_{s=1}^{n+1}\prod_{j=k+1}^{k+s-1}\xi_{j}+\sum_{s=1}^{n/2}\prod_{j=1}^{2s-1}\xi_{k+j}(\xi_{k+2s,k+n}-\xi_{k+2s}),k\ge k_0. \end{align}
  Noticing that $n-2s+1$ is odd, we have $\xi_{k+2s,n+k}-\xi_{k+2s}>0$ and thus the second summation on the right-hand side of the above inequality is positive. We claim that
\begin{align}\label{dps}
  \frac{\sum_{s=1}^{n/2}\prod_{j=1}^{2s-1}\xi_{k+j}(\xi_{k+2s,k+n}-\xi_{k+2s})}{\sum_{s=1}^{n+1}\prod_{j=k+1}^{k+s-1}\xi_{j}} \rightarrow 0, k\ge k_0, \end{align}
uniformly in $k$ as $n\rto.$

For this purpose, as done above, fix $\varepsilon>0$ and let $k_3>0$ be an even number such that $r^{k_3}<\varepsilon.$
Then for $k\ge k_0,$ we have
\begin{align}\label{tsp}
  \sum_{s=1}^{n/2}&\prod_{j=1}^{2s-1}\xi_{k+j}(\xi_{k+2s,k+n}-\xi_{k+2s})\\
  &=\sum_{s=1}^{\frac{n-k_3}{2}}+\sum_{s=\frac{n-k_3}{2}+1}^{\frac{n}{2}}\prod_{j=1}^{2s-1}\xi_{k+j}(\xi_{k+2s,k+n}-\xi_{k+2s})
  =:\mathrm{(III) +(IV)}.\no
\end{align}
Since  $\xi_k>0,k\ge k_0$ and $\lim_{k\rto}\xi_k^{-1}=\varrho\ge1,$  applying  \eqref{uc} and \eqref{up},  we have
\begin{align*}
   \frac{\mathrm{(IV)}}{\sum_{s=1}^{n+1}\prod_{j=k+1}^{k+s-1}\xi_{j}}\rightarrow 0, k\ge k_0
\end{align*}
uniformly in $k$ as $n\rto.$

Now we turn to consider the term (III) on the right-hand side of \eqref{tsp}. It follows from \eqref{up} and Lemma \ref{dxr} that
\begin{align*}
\frac{\mathrm{(III)}}{\sum_{s=1}^{n+1}\prod_{j=k+1}^{k+s-1}\xi_{j}}&\le \sum_{s=1}^{(n-k_3)/2}(\xi_{k+2s,k+n}-\xi_{k+2s})\\
  &\le c\sum_{s=1}^{(n-k_3)/2} r^{n-2s} \le c\sum_{s=k_3/2}^{\infty} r^{2s}\\
  &=cr^{k_3}/(1-r^2)\le c\varepsilon/(1-r^2),k\ge k_0.
\end{align*}
Since $\varepsilon$ is arbitrary, for $k\ge k_0,$ we get ${\mathrm{(III)}}\Big/{\sum_{s=1}^{n+1}\prod_{j=k+1}^{k+s-1}\xi_{j}}\rightarrow 0$ uniformly in $k$ as $n\rto.$
We thus come to the conclusion that \eqref{dps} is true. As a consequence, dividing by $\sum_{s=1}^{n+1} \xi_{k+1}\cdots \xi_{k+s-1}$ on both sides of \eqref{fu} and taking the upper limit, we conclude that
\begin{align}\label{upb}
  \varlimsup_{n\rightarrow\infty}\frac{\sum_{s=1}^{n+1}\prod_{j=k+1}^{k+s-1}\xi_{j,k+n}}{\sum_{s=1}^{n+1}\prod_{j=k+1}^{k+s-1}\xi_{j}}\le 1, k\ge k_0,
\end{align} where the upper limit is taken uniformly in $k.$

For a lower limit, from \eqref{pxi}, \eqref{txi} and \eqref{up}, for $k\ge k_0,$  we get
 \begin{align*}
  \sum_{s=1}^{n+1}\prod_{j=k+1}^{k+s-1}\xi_{j,k+n}&\ge \sum_{s=1}^{n+1}\prod_{j=k+1}^{k+s-1}\xi_{j}+\sum_{s=1}^{n/2}\prod_{j=1}^{2s-2}\xi_{k+j,k+n}(\xi_{k+2s-1,k+n}-\xi_{k+2s-1})\\
  &\ge \sum_{s=1}^{n+1}\prod_{j=k+1}^{k+s-1}\xi_{j}-c\sum_{s=1}^{n/2}\prod_{j=1}^{2s-2}\xi_{k+j}|\xi_{k+2s-1,k+n}-\xi_{k+2s-1}|.
       \end{align*}
  Similarly to \eqref{dps}, we can show that
  \begin{align*}\frac{\sum_{s=1}^{n/2}\prod_{j=1}^{2s-2}\xi_{k+j}|\xi_{k+2s-1,k+n}-\xi_{k+2s-1}|}{\sum_{s=1}^{n+1}\prod_{j=k+1}^{k+s-1}\xi_{j}}\rightarrow 0,k\ge k_0,
      \end{align*} uniformly in $k$ as $n\rto.$
  It thus follows that \begin{align}\label{lowb}
  \varliminf_{n\rightarrow\infty}\frac{\sum_{s=1}^{n+1}\xi_{k+1,k+n}\cdots\xi_{k+s-1,k+n}}{\sum_{s=1}^{n+1}\xi_{k+1}\cdots\xi_{k+s}}\ge 1, k\ge k_0,
\end{align} where the lower limit is taken uniformly in $k.$  Taking \eqref{upb} and \eqref{lowb} together, we see that \eqref{cone} is true whenever $bd>a\theta.$
 Therefore we come to the conclusion that \eqref{cone} is always true so that \eqref{ps} is proved for $i=j=1.$

Finally, we deal with the case $i=1$ and $j=2.$ Notice that
\begin{align}\label{oo}
  \frac{\sum_{s=1}^{n+1}\mb e_1A_{k+s}\cdots A_{k+n}\mb e_2^t}{\sum_{s=1}^{n+1}\xi_{k+s}^{-1}\cdots \xi_{k+n}^{-1}}=\frac{\tilde b_{k+n}\sum_{s=1}^{n}\mb e_1A_{k+s}\cdots A_{k+n-1}\mb e_1^t}{1+\xi_{k+n}^{-1}\sum_{s=1}^{n}\xi_{k+s}^{-1}\cdots \xi_{k+n-1}^{-1}}.
\end{align}
  Since $\lim_{n\rto}\xi_n^{-1}=\varrho\ge1,$ it is easy to see that
  $\xi_{k+n}^{-1}\rightarrow \varrho,$ $\tilde b_{k+n}\rightarrow b$ and $\z(\sum_{s=1}^{n}\xi_{k+s}^{-1}\cdots \xi_{k+n-1}^{-1}\y)^{-1}\rightarrow 0$ uniformly in $k$ as $n\rto.$ Taking the results for $i=j=1$ we proved above, from \eqref{oo}, we conclude that
  $$\frac{\sum_{s=1}^{n+1}\mb e_1A_{k+s}\cdots A_{k+n}\mb e_2^t}{\sum_{s=1}^{n+1}\xi_{k+s}^{-1}\cdots \xi_{k+n}^{-1}}\rightarrow \frac{b}{\varrho}\phi(1,1,k)=\frac{b}{\varrho-\varrho_1}=\phi(1,2,k),k\ge k_0,$$
  uniformly in $k$ as $n\rto.$ We thus finish the proof of \eqref{ps} for $i=1$ and $j=2.$
   \qed
\subsection{Proof of Theorem \ref{pml}}

Theorem \ref{pml} is a direct consequence of Theorem \ref{s}. For \eqref{mx}, we prove here only the case $i=j=1,$ since the other three cases can be proved similarly. By some easy computation, we see that for $k\ge1, 1\le s\le n,$
\begin{align}\label{km}
\mb e_1M_{k+s}\cdots M_{k+n}\mb e_1^t=\mb e_1 A_{k+s}\cdots A_{k+n}\mb e_1^t-\frac{\theta_{k+n+1}}{b_{k+n+1}}\mb e_1A_{k+s}\cdots A_{k+n}\mb e_2^t.
\end{align}
Thus, using Theorem \ref{s} and the facts $\lim_{n\rto}\theta_{n}=\theta$ and $\lim_{n\rto}b_{n}=b,$
we get that uniformly in $k,$
\begin{align*}
 \lim_{n\rto} \frac{\mb e_1M_{k+1}\cdots M_{k+n}\mb e_1^t}{\xi_{k+1}^{-1}\cdots\xi_{k+n}^{-1}}=\frac{\varrho}{\varrho-\varrho_1}-\frac{\theta}{b}\frac{b}{\varrho-\varrho_1}
 =\frac{\varrho-\theta}{\varrho-\varrho_1}, k\ge k_0,
\end{align*}
which finishes the proof of \eqref{mx} for $i=j=1.$

Next we prove \eqref{ms}. If $i=j=1,$ then using \eqref{km} and Theorem \ref{s}, we have
\begin{align*}
  \frac{\sum_{s=1}^{n+1}\mb e_1M_{k+s}\cdots M_{k+n}\mb e_1^t}{\sum_{s=1}^{n+1}\xi_{k+s}^{-1}\cdots \xi_{k+n}^{-1}}&=\frac{\sum_{s=1}^{n+1}\mb e_1 A_{k+s}\cdots A_{k+n}\mb e_1^t}{\sum_{s=1}^{n+1}\xi_{k+s}^{-1}\cdots \xi_{k+n}^{-1}}-\frac{\theta_{k+n+1}}{b_{k+n+1}}\frac{\mb e_1A_{k+s}\cdots A_{k+n}\mb e_2^t}{\sum_{s=1}^{n+1}\xi_{k+s}^{-1}\cdots \xi_{k+n}^{-1}}\\
  &\rightarrow \frac{\varrho}{\varrho-\varrho_1}-\frac{\theta}{b}\frac{b}{\varrho-\varrho_1}
 =\frac{\varrho-\theta}{\varrho-\varrho_1}, k>k_0,
\end{align*}
uniformly in $k$ as $n\rto.$ If $i=1$ and $j=2,$ \eqref{ms} follows trivially from \eqref{ps} since  $\mb e_1M_{k+s}\cdots M_{k+n}\mb e_2^t=\mb e_1A_{k+s}\cdots A_{k+n}\mb e_2^t$ for $k\ge1,$ $1\le s\le n.$ \qed

\vspace{.5cm}

\noindent{{\bf \Large Acknowledgements:}} 
This project is supported by National
Natural Science Foundation of China (Grant No. 11501008; 12071003).


\end{document}